\documentclass[11pt,leqno]{article}
\usepackage{amsmath}
\usepackage{amssymb}
\def\overset#1#2{{\mathrel{\mathop {{#2}_{}}\limits^{#1}}}}
\def\underset#1#2{{\mathrel{\mathop {{}_{} {#2}}\limits_{{#1}_{}}}}}
\def\upplim_#1{\underset{#1}{\overline\lim}\;}
\def\lowlim_#1{\underset{#1}{\underline\lim}\;}
\newtheorem{claim}[equation]{\indent {\it Claim}\rm}

\newtheorem{corollary}[equation]{Corollary}
\newtheorem{definition}[equation]{\indent{\it Definition}\rm }

\newtheorem{lemma}[equation]{Lemma}
\newtheorem{proposition}[equation]{Proposition}
\newtheorem{remark}[equation]{\indent \rm {\it Remark}}

\newtheorem{theorem}[equation]{Theorem}

\newcommand{\Aut}{\mathop{\mathrm{Aut}}}
\newcommand{\Alb}{\mathop{\mathrm{Alb}}}
\newcommand{\C}{{\mathbf{C}}}

\newcommand{\Div}{\mathop{\mathrm{Div}}}
\renewcommand{\div}{\mathop{\mathrm{div}}}

\newcommand{\I}{{\mathcal{I}}}
\newcommand{\id}{\mathrm{id}}

\newcommand{\Lie}{{\mathop{\mathrm{Lie}\,}}}
\newcommand{\Lin}{\mathop{\mathrm{Lin}}}
\newcommand{\m}{\mathfrak{m}}
\newcommand{\Mor}{\mathop{\mathrm{Mor}}}
\newcommand{\mult}{\mathop{\mathrm{mult}}}
\newcommand{\N}{\mathbf{N}}

\renewcommand{\P}{{\mathbf{P}}}
\newcommand{\Pic}{{\mathop{\mathrm{Pic}}}}

\newcommand{\Spec}{{\mathrm{Spec}~}}
\newcommand{\sS}{{\mathcal S}}
\newcommand{\StD}{\mathrm{St}(D)}
\newcommand{\supp}{\mathrm{Supp}\,}
\newcommand{\tensor}{\otimes}
\newcommand{\Z}{\mathbf{Z}}
\newenvironment{proof}{\par{\it Proof.}}{{\it Q.E.D.}\par\vskip3pt}
\numberwithin{equation}{section}

\date{\hfil}
\begin{document}
\baselineskip=18pt
\title{%
Bounds for Curves in Abelian Varieties%
\thanks{  Research supported in part by Grant-in-Aid
for Scientific Research (A)(1), 13304009.}
}
\author{Junjiro Noguchi and J\"org Winkelmann}
\maketitle
\begin{abstract}
A uniform bound of intersection multiplicities of curves and
divisors on abelian varieties is proved by algebraic geometric methods.
It extends and improves a result obtained by A. Buium with a different
method based on Kolchin's differential algebra.
The problem is modeled after the ``$abc$-Conjecture'' of Masser-Oesterl\'e
for abelian varieties over the  function field of a curve.
As an application a finiteness theorem will be proved for
maps from a curve into an abelian variety omitting an ample divisor.
\end{abstract}

\section{Introduction}

Buium [Bu98] proved the following:
\smallskip

\noindent{\bf Theorem A.}
{\it
Let $C$ be a smooth compact curve, let $A$ be  an abelian variety,
and let $D$ be an effective reduced divisor on $A$
which does not contain any translate of a non-trivial
abelian subvariety of $A$.
Then there exists a number $N$ depending
on $C$, $A$ and $D$ such that for every morphism
$f:C\to A$, either $f(C)\subset D$ or the multiplicities
$\mult_x f^*D, x \in C$
of the pulled-back divisor $f^*D$ are uniformly bounded by $N$.
}
\smallskip

This theorem together with his former one [Bu94] was motivated
by the so-called $abc$-Conjecture of Masser-Oesterl\'e on
abelian varieties.
There is also such a bound of multiplicities
in the Second Main Theorem for (transcendental) holomorphic curves
in abelian and semi-abelian varieties by [NWY00] and [NWY99].
In Theorem A, if a finite subset $S$ of $C$ is fixed
and if $f^{-1}D \subset S$ is imposed,
then the bound of such multiplicities immediately follows from
the quasi-projective algebraicity of the corresponding moduli space
of mappings ([N88]).
Therefore, it is a delicate but important point in Theorem A
that $f$ is allowed to take values in $D$ at arbitrary points of $C$.

In this article we provide an algebraic geometric
proof by making use of jet bundles,
a method different from the one used by Buium [Bu98] who relied on
the theory of Kolchin's differential algebra.
By our proof we have the following:
\begin{enumerate}
\item
The condition that the given divisor $D$
contains no translate of non-trivial abelian subvarieties
can be removed, as conjectured in [Bu98]:
\item
The bound $N$ depends only on the numerical data involved;
the genus $g$ of the compact curve $C$, the dimension
of $A$, and the top self-intersection number $D^{\dim A}$.
\end{enumerate}

Moreover, we show that Theorem A permits a generalization
to the  semi-abelian case and we give an application
to a finiteness theorem for
maps from a curve into an abelian variety omitting an ample divisor
(Theorem \ref{appl}).

Our main result is the following:
\smallskip

\noindent
{\bf Main Theorem.}  {\it
Let $A$ be a semi-abelian variety
(i.e.,
a connected commutative reductive algebraic group),
let $A\hookrightarrow\bar A$ be
a smooth equivariant algebraic compactification,
let $\bar D$ be an effective reduced ample divisor on $\bar A$,
let $D=\bar D\cap A$, and let $C$ be a smooth algebraic curve
with smooth compactification $C\hookrightarrow\bar C$.
Then there exists a number $N\in\N$ such that for every
morphism $f:C\to A$ either $f(C)\subset D$ or
$\mult_x f^*D \le N$ for all $x\in C$.
Furthermore, the number $N$ depends only on the numerical data involved
as follows:
\begin{enumerate}
\item
The genus of $\bar C$ and the number
$\#(\bar C\setminus C)$ of the boundary points of $C$,
\item
the dimension of $A$,
\item
the toric variety (or, equivalently, the associated ``fan'')
which occurs as closure of the orbit in $\bar A$ of the maximal
connected linear algebraic
subgroup $T \cong (\C^*)^t$ of $A$,
\item
all intersection numbers of the form
$D^h\cdot B_{i_1}\cdots B_{i_k}$,
where the $B_{i_j}$ are closures of $A$-orbits in $\bar A$
of dimension $n_j$ and $h+\sum_j n_j=\dim A$.
\end{enumerate}
}
\smallskip

In particular, if we let $A$, $\bar A$, $C$ and $D$ vary within a
flat connected family, then we can find a uniform bound for $N$.
For abelian varieties this specializes to the following result:

\begin{theorem}
\label{thm1.1}
There is a function  $N:\N\times\N\times \N\to\N$
such that the following statement holds:
Let $C$ be a smooth compact curve of genus $g$,
let $A$ be an abelian variety of dimension $n$,
let $D$ be an ample effective divisor on $A$ with
intersection number $D^n=d$.
Let $f:C\to A$ be a morphism.
Then either $f(C)\subset D$ or $\mult_x f^*D\leq N(g,n,d)$
for all $x \in C$.
\end{theorem}

{\em Remark.}
If $D$ is an effective divisor on a compact complex torus $M$ which is
not ample, then there exists a subtorus $\StD$ of $A$ stabilizing $D$
such that $D$ is the pull-back of an ample divisor on $A/\StD$.
Thus one can easily generalize our theorem to the non-ample case.

For toric varieties we obtain
\begin{theorem}
Let $\bar A$ be a toric variety with open orbit $A$,
let $C$ be an affine curve, let $D$ be an ample divisor on $\bar A$.
Then there exists a number $N$ depending on the genus of
the smooth compactification $\bar C$ of $C$, $\#(\bar C\setminus C)$
and the numerical data of $(A,\bar A,D)$ such that for every
morphism $f:C\to A$
we have either $f(C) \subset D$ or $\mult_x f^*D \leq N, x \in C$.
\end{theorem}

{\it Acknowledgement.}  The authors are grateful to Professor
Yukitaka Abe for the discussion on Theorem 
\ref{three-very-ample}; he gave another
proof for the theorem when the algebraic torus $\C^t$ is
compactified by $(\P^1(\C))^t$.

\section{Basic idea}

The goal is to achieve a bound on the possible multiplicities for $f^*D$
where $D$ is a divisor on an abelian variety $A$,
$C$ is a curve and $f$ runs through all morphisms $f:C\to A$
except those with $f(C)\subset \supp D$ (= the support of $D$).

A rather naive idea would be to define subsets $Z_k$ of the space
$\Mor(C,A)$ of all morphisms from $C$ to $A$
by requiring that $f\in Z_k$
if and only if $f^*D$ has multiplicity $\ge k$ somewhere.
Then the $Z_k$ form a decreasing sequence of subsets of $\Mor(C,A)$.
Now, if $\Mor(C,A)$ were an algebraic variety and $Z_k$ were algebraic
subvarieties, this sequence $Z_k$ would eventually have to stabilize,
i.e., there would be a number $N$ such that $Z_k=Z_N$ for all $k\ge N$.
Unfortunately $\Mor(C,A)$ is not necessarily an algebraic variety,
but it may have infinitely many irreducible components.
This infinitude make things complicated.

Our approach is to embed $\Mor(C,A)$ into a larger connected space
$H$,
such that the $Z_k$ can be extended to a sequence of algebraic subsets
$\tilde Z_k$ which eventually stabilizes.
This is done in the following way:
Every morphism $f:C\to A$ induces a morphism from the
Albanese $\Alb(C)$ to $A$ which in turn induces a holomorphic
map between the respective universal coverings.
The crucial fact here is that a morphism between compact complex
tori necessarily lifts to an {\em affine-linear} map
between the corresponding universal coverings.
Clearly the space $\Lin(\C^n,\C^m)$ of all affine-linear maps from
one $\C^n$ to a $\C^m$ is connected,
and also carries a natural structure as an algebraic variety.

Now let $H=C\times D\times \Lin(\C^n,\C^m)$.
Then we can define subsets $\tilde Z_k\subset H$ in the following way:
$(p,q,\phi)\in\tilde Z_k$ if and only if
there is a holomorphic map germ
$f:(C,p)\to(A,q)$ such that $f^*D$ has multiplicity $\ge k$ at $p$
and $f$ lifts to an affine-linear  map
whose linear part is $\phi$.
Then $H$ is an algebraic variety, and the subsets $\tilde Z_k$
form a descending sequence of closed algebraic subvarieties,
which eventually stabilizes.

Note that, up to the choice of a base point, the universal
covering of the Albanese torus of a projective variety
can be canonically identified with the dual vector space of
the space of global holomorphic $1$-forms.
Correspondingly, in the proof itself we discuss induced
maps between these vector spaces rather than maps between
the universal covering of the respective Albanese tori.

\section{Basic notions}

\subsection{Commutative reductive groups}

We recall some basic terminologies.
\begin{definition}{\rm
\begin{enumerate}
\item
A complex Lie group $G$ is said to be \begin{em}reductive\end{em}
if $G$ itself is the only complex Lie subgroup $H\subset G$ containing a
maximal compact subgroup of $G$.
\item
A (complex) \begin{em} semi-torus\end{em} is a
connected commutative reductive complex Lie group.
\end{enumerate}}
\end{definition}

Every semi-torus $A$ admits a short exact sequence of
commutative complex Lie groups
\begin{equation}
\label{pres}
1 \to T \to A \to M \to 1,
\end{equation}
where $M$ is a compact complex torus and $T\cong \left(\C^*\right)^t$
for some $t \in \N$.

A commutative algebraic group $T$ is called
a \begin{em} semi-abelian variety\end{em}
if there exists a short exact sequence of algebraic groups
as in (\ref{pres}) (with $M$ being an abelian variety).

Neither presentation (\ref{pres}),
nor the compact torus $M$ is unique in the complex-analytic category.
However, if $A$ is a semi-abelian variety,
then there is one and only one such algebraic presentation.
This fact can be deduced in the following way:
The theorem of Chevalley implies that
there is a linear algebraic subgroup $T$ of $A$
such that $M=A/T$ is an abelian variety ([Bo91]).
This linear algebraic subgroup $T$ is unique, because every
morphism from a linear algebraic group to an abelian variety
is constant.
Thus, for a semi-abelian variety $A$ there is exactly one
presentation
\[
1 \to T \to A \to M \to 1,
\]
such that $T$, $M$ and the morphisms are all algebraic.

A different way to characterize semi-tori is the following:
A complex Lie group $A$ is a semi-torus if and only if it is isomorphic
to a quotient $(\C^n,+)/\Gamma$, where $\Gamma$ is a discrete
subgroup of $\C^n$ spanning $\C^n$ as complex vector space.

\begin{definition}\label{def-qa}{\rm
Let $A$ be a complex semi-torus.
An equivariant smooth compactification
$A\hookrightarrow \bar A$
is said to be \begin{em}quasi-alge\-bra\-ic\end{em}
if
every isotropy subgroup of the $A$-action on $\bar A$ is reductive.}
\end{definition}

If $A$ is a semi-abelian variety, then
evidently every algebraic compactification $\bar A$ is quasi-algebraic
(but not conversely).

\begin{lemma}\label{snc}
Let $A\hookrightarrow\bar A$ be a smooth quasi-algebraic compactification
of a complex semi-torus and let $D$ be an (effective, reduced)
$A$-invariant divisor on $\bar A$.
Then $D$ is a divisor with only simple normal crossings.
\end{lemma}
\begin{proof}
Let $p\in \supp D$ and $A_p=\{a\in A:a(p)=p\}$.
By assumption, $A_p$ is reductive. This implies that the $A_p$-action
can be linearized near $p$,  i.e. there is an $A_p$-equivariant
biholomorphism between an open neighborhood of $p$
in $\bar A$ and an open neighborhood of $0$ in the vector space
$V=T_p\bar A$.
It follows that a neighborhood of $p$ in $D$ is isomorphic
to a union of vector subspaces of codimension one in $V$. From
the assumption that the $A_p$-action is effective, one can
deduce that this is a transversal union, i.e.~$D$ is a simple normal
crossings divisor near $p$.
\end{proof}

\subsection{Toric varieties}\label{toric}

An equivariant compactification of $T=(\C^*)^t$ is called a
\begin{em}toric variety\end{em}.
Here we are interested only in toric varieties which are
nonsingular, compact and moreover projective.
We recall some well-known properties of such toric varieties
(see, e.g., [O85]).

Let $\bar T$ be a projective smooth equivariant compactification
of $T=(\C^*)^t$.
Then $\bar T$ is simply connected and the homology
$H_*(\bar T, \Z)$ is generated (as $\Z$-module) by the
fundamental classes of closures of $T$-orbits.
In particular, the Chern classes of the tangent bundle
are linear combinations
of fundamental classes of orbit closures.

There are only finitely many $T$-orbits, and each $T$-invariant
divisor has only simple normal crossings.
Every ample line bundle on $\bar T$ is already very ample.
Note that
$H^1(\bar T ,{\mathcal O}_{\bar T})=0$ and hence that
every topologically trivial line bundle on $\bar T$
is already holomorphically trivial.

There are only countably many non-isomorphic
compact toric varieties.
Toric varieties can be classified by certain combinatorial
data, called {\em fans} (cf.\ [O85]).

\begin{lemma}\label{neu}
Let $\bar T$ be a smooth projective toric variety.
Then there exists a number $k\in\N$ such that for
every ample line bundle $H$
and every reduced $T$-invariant hypersurface $D$
the line bundle $H^{k}\tensor L(-D)$ is ample.
\end{lemma}
\begin{proof}
There are only finitely many $T$-invariant reduced
hypersurfaces $D_1,\ldots,D_r$ on $\bar T$
and only finitely many $T$-invariant curves
$C_1,\ldots,C_s$.
We choose $k$ such that $k>D_i\cdot C_j$
for all $i\in\{1,\ldots,r\}$, $j\in\{1,\ldots,s\}$.
Then $\left(k c_1(H)-D_i\right)\cdot C_j>0$ for all $i,j$.
By the toric Nakai criterion (see [O85], Theorem 2.18) it follows that
$H^{k}\tensor L(-D_i)$ is ample for every $i$.
\end{proof}

\subsection{Compactifications of semi-tori}
Let $A$ be a semi-torus with presentation
\begin{equation}\label{presentation}
1 \to T \to A \to M \to 1,
\end{equation}
where $M$ is a compact complex torus and $T\cong(\C^*)^t$.
An equivariant compactification $\bar A\hookleftarrow A$ is called
{\em compatible} with the given presentation, if the
projection $\pi:A\to M$ extends to a holomorphic map
$\bar\pi:\bar A\to M$.
In this case, $\bar T=\bar\pi^{-1}(1)$
 is a \begin{em}toric variety\end{em}.
Furthermore
\[
\bar A= A\times_T \bar T.
\]
(Here $A\times_T \bar T$ denotes the quotient of $A\times\bar T$
under the diagonal action of $T$.)
We say that $A\hookrightarrow\bar A$ is a compactification
of type $\bar T$.

Let us now discuss the case where $A$ is a semi-abelian varierty.
We claim that in this case every smooth algebraic equivariant
compactification $A\hookrightarrow\bar A$ is compatible
with the unique algebraic presentation.
Indeed, let $A\hookrightarrow\bar A$ be a smooth equivariant
algebraic compactification.
Since $A$ acts effectively, for every $x\in\bar A$ and
$g\in A_x=\{g\in A; g(x)=x\}\setminus\{1\}$ the automorphism of $\bar A$
given by $g$ induces a non-trivial automorphism of
the local ring  at $x$ and therefore on the space of
$k$-jets (see the next subsection 3.4) at $x$ for some $k\in\N$.
It follows that $A_x$ admits a faithful linear representation.
Thus $(A_x)^0\subset T$.
Therefore every algebraic equivariant compactification is
compatible with the unique algebraic presentation.

\begin{lemma}\label{trans-U}
Let $A\hookrightarrow\bar A$ be a smooth equivariant compactification
compatible with a presentation as in (\ref{presentation}).
Let $U$ be the maximal compact subgroup of $T$.
Then $\bar A\to M$ can be described as a locally trivial fiber bundle
whose transition functions are locally constant functions with values in
$U$.
\end{lemma}
\begin{proof}
Let $K$ be a maximal compact subgroup of $A$, and $V=\Lie K\cap i\Lie K$.
For a sufficiently small open neighborhood $U$ of $0$ in $V$ the
map $U\ni v\mapsto\exp(v)\cdot e_M$ is a biholomorphism onto an open
neighborhood of the unit $e_M \in M$.
Every $q\in M$ can be written in the form
$q=\pi(k)$ for some $k\in K$.
Then $U\ni v\mapsto\exp(v)\cdot \pi(k)$ parameterizes an open
neighborhood of $q$ in $M$ and the diagram
\[
\begin{matrix}
(U\times\bar T)\ni(v,x)& \mapsto & \exp(v)\cdot k \cdot x \\
\downarrow {}_{\mathrm{pr}_1} && \downarrow \\
U \ni v & \mapsto & \exp(v) \cdot \pi(k) \\
\end{matrix}
\]
gives a local trivialization for this neighborhood.
Everything is canonical except for the choice of $k$.
Since $\pi(k)$ must equal $q$, the element $k$ is well-defined
up to multiplication by an element in $U=K\cap\pi^{-1}(e_M)$.
Hence the local trivializations constructed in this way
realize $\bar\pi$ as a fiber bundle whose transition functions
are locally constant with values in $U$.
\end{proof}

\subsection{Sheaf of logarithmic 1-forms}

Let $N$ be a complex manifold, and let $D$ be a divisor on $N$
with only normal crossings.
For $p \in N$ we take a local coordinate neighborhood
$U(x_1, \ldots, x_n)$ such that
$D \cap U=\{x_1\cdots x_k=0\}$ ($0 \leq k \leq n$).
Then the sheaf $\Omega^1_N(\log D)$ of logarithmic $1$-forms
along $D$ is the sheaf of those meromorphic $1$-forms
 locally generated over the structure
sheaf $\mathcal{O}=\mathcal{O}_N$ by
$$
\frac{dx_1}{x_1}, \ldots, \frac{dx_k}{x_k}, dx_{k+1},
\ldots, dx_n.
$$
The sheaf $\Omega^1_N(\log D)$ is locally free, and the
dual bundle of its associated vector bundle is called
the logarithmic tangent
bundle, denoted by $\mathbf{T}(N, \log D)$ ([I76]).

Let $M,N$ be complex manifolds and $D,E$ divisors with
only normal crossings on $M$ resp.~$N$.
Let $F:M\setminus D\to N\setminus E$ be a holomorphic
map which extends to a meromorphic map from $M$ to $N$.
Then $F^*\Omega^1_N(\log E)\subset \Omega^1_M(\log D)$.

In general, let $X$ be a reduced complex space and let $E$
be a reduced complex subspace of $X$
(we will actually deal only with algebraic $X$).
By the well-known Hironaka resolution there is a proper
holomorphic mapping $\lambda: \tilde X \to X$ such that
$\tilde X$ is smooth, and $\tilde E=\pi^{-1}(E)$ is
a divisor with only normal crossings.
Then we have $\Omega^1_{\tilde X}(\log \tilde E)$
defined as above.

We define the sheaf of logarithmic 1-forms on $X$ along $E$ by
its direct image sheaf:
\begin{equation}
\label{log-form}
\Omega^1_{X}(\log E)=:R^0 \lambda_* \Omega^1_{\tilde X}(\log \tilde E).
\end{equation}
Let $\lambda': \tilde X' \to X$ be another resolution such that
$\tilde X'$ is smooth and $\tilde E'={\lambda'}^{-1}E$ is
a divisor with only normal crossings.
Then
$$
R^0 \lambda_* \Omega^1_{\tilde X}(\log \tilde E)\cong
R^0 \lambda'_* \Omega^1_{\tilde X'}(\log \tilde E').
$$
(see [I76]).
Hence $\Omega^1_{X}(\log E)$ is well-defined.

Let $X$ be a (possibly non-compact)
smooth complex algebraic variety.
By [I77] we have the {\it quasi-Albanese variety} $\Alb(X)$
and the {\it quasi-Albanese mapping} $\alpha_X: X \to \Alb(X)$.
If $X$ is compact, the quasi-Albanese variety is just the usual
Albanese variety. In general, the quasi-Albanese variety $\Alb(X)$
is a semi-abelian variety.
As in the case of Albanese varieties, we have a universal
property: For a semi-Abelian variety $A$ and a morphism
$f: X \to A$, there is a unique morphism $\tilde f : \Alb(X) \to A$
such that $f=\tilde f \circ \alpha_X$.

\subsection{Jet bundles}

Let $\Delta=\{|z|<1\} \subset \C$ be the unit disk with center at the
origin.
The space $J^k_p(X)$ of $k$-jets at a point $p$ of a complex manifold $X$
can be defined as the quotient space of all holomorphic mappings
of the space germ $(\Delta, 0)$ into the space germ $(X, p)$
by the equivalence relation given as follows:
$f \overset{k}{\sim} g$ if their Taylor series expansions
agree up to degree $k$.
Set
$$
J^k(X)=\bigcup_{p\in X} J^k_p(X),
$$
which is called the $k$-jet bundle over $X$.
For $k=1$ there is a natural isomorphism $J^1(X)\cong {\bf T}(X)$
with the holomorphic tangent bundle ${\bf T}(X)$ over $X$.
For $k \geq 2$ $J^k(X)$ carries no natural vector bundle structure.
However, the $\C^*$-action on $(\Delta,0)$
induces a natural $\C^*$-action on $\Spec\C\{t\}/(t^{k+1})$ via
$$
c \cdot j(t)=j(c t), \quad c \in \C^*,~j(t)\in
\Spec\C\{t\}/(t^{k+1}),
$$
and hence on $J^k_p(X)$.
This $\C^*$-action can be used to define the notion of polynomials
of weighted degree: A function $P$ on $J^k_p(X)$ is called a
{\em polynomial of weighted degree $d$}
if $P( c \cdot j)=c^dP(j)$.
The points of $J^k_p(X)$ are separated by polynomials of
weighted degree $\le k$.

Let $Y \subset X$ be an analytic subspace with defining ideal
sheaf $\I_Y$.
For $p \in Y$ we define the $k$-jet space of $Y$ by
$$
J^k_p(Y)=\{j \in J^k_p(X); (j^*\I_Y)_0 \subset (t^{k+1})\},\quad
J^k(Y)=\bigcup_{p\in Y}J^k_p(Y).
$$
Then $J^k(Y)$ is a complex subspace of $J^k(X)$.
When $Y$ is non-singular, $J^k(Y)$ coincides the one defined
as above with $X=Y$.
Because $J^k(Y)$ is independent of the used embedding
$Y\hookrightarrow X$, the $k$-jet space $J^k(Y)$ over $Y$
is defined for a (reduced) complex space $Y$.

Equivalently, for any complex space $X$
we can define the space $J^k_p(X)$  of $k$-jets at a point $p\in X$
as the space of holomorphic
mappings from $\Spec \C\{t\}/(t^{k+1})$ to $X$ which map the
geometric point of $\Spec\C\{t\}/(t^{k+1})$ to $p$.

The maximal ideal of $\Spec\C\{t\}/(t^{k+1})$ is denoted by
$\m(k)$.

A holomorphic mapping $\phi:X\to Y$ between complex spaces induces
$\C^*$-equivariant mappings
$$
d^k \phi:J^k(X)\to J^k(Y), \quad
\hbox{i.e., } c d^k \phi(j)=d^k\phi(c j)
$$
for $j\in J^k_p(X)$, $c \in\C^*$.

Let $A$ be a semi-abelian variety, and
let $\bar A$ be a projective
algebraic equivariant compactification.
Then $A\hookrightarrow\bar A$ is
quasi-algebraic in the sense of Definition \ref{def-qa}.
Hence the $A$-invariant divisor $\partial A$ has only simple
normal crossings (Lemma~\ref{snc}).
We have the logarithmic tangent bundle
$\mathbf{T}(\bar A, \log\partial A)$
and logarithmic jet bundles
$J^k(\bar A, \log\partial A)$ over $\bar A$ ([I76], [N86]).
For $k=1$ we have
\begin{equation}
\label{k=1}
J^1(\bar A, \log\partial A)\cong \mathbf{T}(\bar A, \log\partial A).
\end{equation}
Moreover, there is a natural bundle morphism,
\begin{equation}
\label{jet}
J^k(\bar A, \log\partial A) \to J^k(\bar A).
\end{equation}

\begin{proposition}
\label{jet-trivial}
Let $A\hookrightarrow\bar A$ be a quasi-algebraic smooth equivariant
compactification of a complex semi-torus $A$.
Then the logarithmic jet bundles $J^k(\bar A, \log\partial A)$
are trivialized as
$$
J^k(\bar A, \log\partial A) \cong \bar A \times(\m(k)\tensor\Lie A).
$$
\end{proposition}

\begin{proof}
We first show that the logarithmic tangent bundle is trivialized by the
$A$-fundamental vector fields.
Since we required the compactification to be
quasi-algebraic in the sense of Definition \ref{def-qa},
all the isotropy groups are reductive.
Therefore
the action of every isotropy subgroup is linearizable in some
open neighborhood.
Thus for every point on $\bar A$ we can find
a system of local coordinates in which the $A$-fundamental
vector fields are simply
\[
\frac{\partial}{\partial z_1},\ldots,\frac{\partial}{\partial z_h},
z_{h+1}\frac{\partial}{\partial z_{h+1}},\ldots,
z_{n}\frac{\partial}{\partial z_n},
\]
for some $h,n\in\N$,
where $n-h$ equals the dimension of the isotropy group.
Hence we have
\begin{equation}
\label{tan-lie}
{\bf T}(\bar A, \log\partial A)\cong \bar A \times \Lie A.
\end{equation}

For $k\geq 2$, (\ref{tan-lie}) induces a global trivialization of
$J^k(\bar A, \log\partial A)$ in the following way:
We regard the $k$-jets over $A$ as maps from $S_k=\Spec\C\{t\}/(t^{k+1})$
to $A$. As before, let $\m(k)$ denote the maximal ideal $(t)$
of $\C\{t\}/(t^{k+1})$.
For $p\in A$ and
$\alpha=\sum_i \alpha_i\tensor v_i\in \m(k)\tensor\Lie A$
we define a map (germ) from $\Spec\C\{t\}/(t^{k+1})$
by
\[
\alpha: t \mapsto \exp\left( \sum_i \alpha_i(t)v_i\right)\cdot p .
\]
A calculation in local coordinates shows that this gives
a trivialization,
$J^k(\bar A, \log\partial A)|_A \cong A \times(\m(k)\tensor\Lie A)$.
This trivialization holomorphically extends over $\bar A$
(see [N86]).
\end{proof}

Now fix a point $p \in \bar A$ and consider the induced
map $\m(k) \tensor\Lie A\to J_p^k(\bar A, \log\partial A)$.
Let $V$ be a complex vector space and let $E$ be
a vector space of linear mappings from $V$ to $\Lie A$.
Then we obtain a map
$\m(k) \tensor V\times E\to J^k_p(\bar A,\log\partial A)$
induced by the natural pairing $V\times E \to\Lie A$.
Observe that this map
$\m(k) \tensor V\times E \to J^k_p(\bar A,\log\partial A)$
is polynomial of degree $\le k$ in $E$.

\section{Proof of the Main Theorem in absolute case}

Here we deduce the following key lemma.

\begin{lemma}
\label{key}
Let $A$, $B$ be complex semi-abelian varieties
with smooth equivariant algebraic
compactifications $\bar A$ and $\bar B$.
Let $\phi: B \to A$ be a morphism.
Then there is a linear map $\lambda_\phi \in \Lin(\Lie B,\Lie A)$
such that in view of Proposition \ref{jet-trivial}
\begin{equation}
\begin{array}{rccc}
d^k\phi:& J^k(\bar B,\log\partial B)|_B &
{\to} & J^k(\bar A, \log\partial A)|_A \\
 & \| & \circlearrowleft & \| \\
(\phi, \mathrm{id_{\m(k)}}\tensor\lambda_\phi):&
B \times(\m(k)\tensor\Lie B) &
{\to} & A \times(\m(k)\tensor\Lie A).
\end{array}
\end{equation}
\end{lemma}

\begin{proof}
The morphism $\phi$ induces a bundle morphism
$\phi_*:{\bf T}(\bar B, \log \partial B)\to
{\bf T}(\bar A, \log \partial A)$.
It follows from (\ref{k=1}) and Proposition \ref{jet-trivial}
that there is a Lie algebra homomorphism $\lambda_\phi:\Lie B \to\Lie A$
with $d \phi =(\phi, \mathrm{id_{\m(1)}}\tensor\lambda_\phi)$.
Since the Lie algebras are commutative, a Lie algebra
homomorphism is simply a linear map.
The case of $k\geq 2$ follows from this and Proposition \ref{jet-trivial}.
\end{proof}

Let the notation be as in the Main Theorem.
Let $B$ be the quasi-Albanese variety of $C$
and $\alpha:C \to B$ be the quasi-Albanese mapping.
For the Main Theorem, it suffices to deal with the case when
$\dim B>0$. 

One may assume that $C$ is embedded into $B$, and hence a
(locally closed) subvariety
of $B$.
Let $H=\Lin(\Lie B, \Lie A)$.
We define
$Z_k$ as the set of all $(\tilde\phi ,p,q)$ in $H\times C \times D$
such that
\[
d^k\phi
 (J_p^k( C))
\subset J_q^k( D)
\]
for every holomorphic map-germ $\phi:(C,p)\to (A,q)$
with the property that $d\phi(v)=\tilde\phi(v)$ near $p$ for every
element $v\in\Lie B$, regarded as vector field on $C$.

We claim that the $Z_k$ form a descending sequence of closed
algebraic subsets.
Indeed, $J^k(C)$ (resp.\ $J^k(D)$) is a closed algebraic subvariety
of $J^k(B)|_C\cong C\times \left(\m(k)\tensor\Lie B\right)$
(resp.\ $J_k(A)$).
Thus $J^k_p(C)$ (resp.\ $J^k_q(D)$) can be regarded
as a closed subvariety of $\m(k)\tensor\Lie B$
(resp.\ $\m(k)\tensor\Lie A$).
Then
\begin{equation}\label{k-locus}
Z_k=\{(\lambda,p,q)\in H\times C\times D :
\left(\mathrm{id}_{\m(k)}\tensor\lambda\right)
\left(J^k_p(C) \right)\subset J^k_q(D)\},
\end{equation}
and the algebraicity of the sets $Z_k$ becomes evident.

By Noetherianity of the algebraic Zariski topology
we may conclude that there is a number $N$
such that $Z_k=Z_N$ for all $k\ge N$.

For any morphism $\phi:C\to A$ let $\tilde\phi$ denote the associated
linear map in $H$.
If there exists a pair $(p,q)\in C\times D$ such that $\phi(p)=q$
and $(\tilde\phi,p,q)\in Z_N$, then $(\tilde\phi,p,q)\in Z_k$ for all
$k\ge N$ and consequently $\mult_p \phi^*D\ge k$ for all $k$ which in
turn implies that $\phi(C)\subset D$ by identity principle.
On the other hand, if there is no such pair,
then $\mult_p \phi^*D<N$ for all  $p\in C$.

This finished the proof for the Main Theorem except for the
dependence of the number $N$.
In the rest of this paper, we will prove the said numerical
dependence of $N$ by constructing certain parameter spaces of the
objects which appeared in the above proof.

\section{Line bundles on semi-abelian variety}
\subsection{Notation}
Throughout this section $A$ denotes a semi-abelian variety.
There is a short exact sequence of morphisms of
algebraic groups
\[
1 \to T \to A \stackrel{\pi}{\to} M \to 1,
\]
where $T\cong \left( \C^*\right)^t$, $t\in\N\cup\{0\}$
and $M$ is an abelian variety.
Let $m=\dim M$ and $n=\dim A$.
Let $T\hookrightarrow\bar T$ denote a smooth projective
equivariant algebraic
compactification of $T\cong \left( \C^*\right)^t$; i.e.~$\bar T$ is a
toric variety.
Let $A\hookrightarrow\bar A$ denote a smooth algebraic equivariant
compactification of type $\bar T$,
compatible with the above presentation;
i.e.~$\pi:A\to M$ extends to
an equivariant holomorphic map $\bar\pi:\bar A\to M$.

\subsection{Line bundles}
The $A$-orbits in $\bar A$ of a given codimension
are in one-to-one correspondence
with the $T$-orbits of that codimension in $\bar T$.
In particular, there is a bijective correspondence between
$T$-invariant divisors on $\bar T$ and $A$-invariant divisors
on $\bar A$.
Since $A$ acts with only finitely many orbits on $\bar A$,
a divisor on $\bar A$ is $A$-invariant
if and only if its support is contained in the boundary
$\partial A$.

The presentation
\[
1 \to T \to A\, \overset{\pi}{\to}\, M \to 1
\]
realizes $A$ as a $T$-principal bundle over $M$.
Choosing a vector subspace $W$ of $\Lie A$ transversal to $\Lie T$,
we obtain a splitting of the holomorphic tangent bundle
$\mathbf{T}(\bar A)=H_{\bar A}\oplus V_{\bar A}$
into horizontal and vertical tangent bundles,
where $H_{\bar A}$ is the subbundle spanned by the fundamental
vector fields coming from $W$ and $V_{\bar A}$ is given as the
kernel of $d \bar\pi:T_{\bar A}\to T_M$.
Note that $H_{\bar A}$ is a bundle with
flat connection, because $W$ is a commutative subalgebra of the Lie
algebra $\Lie A$.

We choose a K\"ahler form $\Omega_0$ on $\bar T$.
By averaging, we may assume that $\Omega_0$ is $U$-invariant 
where $U$ denotes
the maximal compact subgroup of $A$.
Recall that $\bar\pi:\bar A\to M$ can be realized as a fiber bundle
with $U$ as structure group (Lemma \ref{trans-U}).
Using these facts it is clear that
$\Omega_0$ induces a hermitian metric on the bundle $V_{\bar A}$
of vertical tangent vectors.
Extending it by zero on $H_{\bar A}$ we obtain a closed
semi-positive $(1,1)$-form $\Omega$ on $\bar A$
which is K\"ahler on each fiber of
$\bar\pi:\bar A\to M$.
In this way we have

\begin{lemma}
\label{semikaehler}
There is a closed semi-positive $(1,1)$-form $\Omega$ on $\bar A$ such
that
\begin{enumerate}
\item
$\Omega$ is positive definite on every fiber of $\bar A\to M$,
\item
if $Z$ is an $A$-invariant divisor, then $\Omega^t|_Z\equiv 0$
(where $t=\dim T$).
\end{enumerate}
\end{lemma}

For a divisor $E$ on $\bar A$ we denote by $L(E)$ its
associated line bundle on $\bar A$.

\begin{lemma}
\label{eff-div}
Let $Z$ be a divisor on $\bar A$ such that $Z\cap A$ is effective.
Assume that there are a line bundle $L_0 \in \Pic (M)$
and a divisor $E$ with $\supp E \subset \partial A$,
satisfying
$L(Z)\tensor \bar\pi^*L_0^{-1} \cong L(E)$.
Then $c_1(L_0)\ge 0$.
\end{lemma}

\begin{proof}
Assume the contrary.
Recall that $M$ is a compact complex torus with universal covering
$\pi_0:\C^m \to M$ (with $m=\dim M$).
We may regard the Chern class $c_1(L_0)$ as
bilinear form on the vector space $\C^m$.
Since $c_1(L_0)$ is not semi-positive definite,
there is a vector $v\in\C^m$ with $c_1(L_0)(v,v)<0$.
Set $W=\{w \in \C^m:c_1(L_0)(v,w)=0\}$.
Let $\mu$ be a semi-positive $(1,1)$-form on $\C^m$ such that
$\mu(v,\cdot)\equiv 0$ and $\mu|_{W\times W}>0$.
Let $\Omega$ be as in Lemma \ref{semikaehler},
and consider the $(n-1,n-1)$-form $\omega$ on $\bar A$ given
by
\begin{equation}
\label{4.3}
\omega=\Omega^t \wedge\bar \pi^*\mu^{m-1}, \quad m=\dim M.
\end{equation}
By construction we have $\omega\wedge\bar\pi^*c_1(L_0) <0$.

Let $Z=Z'+Z''$ so that $Z'$ is effective and no
component of $Z'$ is contained in $\partial A$,
and $\supp Z'' \subset \partial A$.
By the Poincar\'e duality,
$$
\int_{\bar A}c_1(L(Z))\wedge\omega = \int_{Z} \omega.
$$
Since $\omega\wedge c_1(L(E))=0$, we  have
$$
\int_{\bar A}c_1(L(Z))\wedge\omega=
\int_{\bar A}\bar\pi^*c_1(L_0)\wedge\omega <0.
$$
On the other hand,
$$
\int_{Z}\omega = \int_{Z'}\omega + \int_{Z''} \omega.
$$
Note that
$\int_{Z'}\omega\ge 0$,
because $Z'$ is effective and $\omega\ge 0$,
and that $\int_{Z''} \omega=0$, because $\supp Z'' \subset \partial A$,
and $\Omega^t$ vanishes on $\partial A$ by Lemma \ref{semikaehler} (ii).
Thus we deduced a contradiction.
\end{proof}

\begin{proposition}
\label{decomp}
\begin{enumerate}
\item
For every line bundle $L$ on $\bar A$ there exists a line bundle
$L_0$ on $M$ and an $A$-invariant divisor $D$ on $\bar A$ such that
$L=\bar\pi^* L_0\tensor L(D)$.
\item
Let $L_1, L_2\in\Pic(M)$ and
let $D_1,D_2$ be $A$-invariant divisors on $\bar A$.
Assume that
\[
\bar\pi^* L_1\tensor L(D_1) \cong \bar\pi^* L_2\tensor L(D_2).
\]
Then there is a topologically trivial line bundle $H\in\Pic^0(M)$
such that $L_1=L_2\tensor H$ and $\pi^*H\cong L(D_2-D_1)$.
\end{enumerate}
\end{proposition}

\begin{proof}
(i)  Recall that for the toric variety $\bar T$ the cohomology group
$H^2(\bar T,\Z)$ is generated by the fundamental classes of
invariant divisors. Therefore there exists an $A$-invariant divisor
$D$ on $\bar A$ such that $L\tensor L(D)$ is topologically
trivial on $\bar T$.
For a line bundle on a toric variety, the topological
triviality implies the holomorphic triviality.
Hence, $L \tensor L(D)$ is trivial along the $\bar\pi$-fibers.
Therefor $R^0\bar \pi_*{\mathcal O}(L\tensor L(D))$
is an invertible coherent sheaf on $M$ and thus
$L=\bar\pi^*L_0\tensor L(D)$ for some line bundle $L_0 \in \Pic(M)$.

(ii)  Since $(D_1-D_2)\cap A=\emptyset$ and
$L(D_1-D_2)\cong \bar \pi^*(L_2\tensor L_1^{-1})$,
we may apply Lemma~\ref{eff-div}
to $Z=D_1-D_2$ and $L_0=L_2\tensor L_1^{-1}$.
Hence $c_1(L_2\tensor L_1^{-1})\ge 0$.
Similarly, $c_1(L_1\tensor L_2^{-1})\ge 0$.
Thus $c_1(L_2\tensor L_1^{-1})=0$.
Because $H^2(M,\Z)$ is torsion-free ($M$ is a torus),
$L_2\tensor L_1^{-1}$ is topologically trivial.
By the some arguments as in (i), this implies that this
bundle comes from $M$.
\end{proof}

\begin{proposition}
\label{ample-decomp}
Let $L_0$ be a line bundle
on $M$ and let $D$ be an $A$-invariant divisor on $\bar A$.
Let $L=\bar\pi^* L_0\tensor L(D)$.
Then $L$ is ample on $\bar A$
if and only if both $L_0$ and $D\cap \bar T$ are ample.
\end{proposition}

\begin{proof}
The ampleness of $L_0$ follows from that of
$L$ in the same way as
in the proof of
Lemma~\ref{eff-div}.
Furthermore, $L|_F\cong L(D)|_{\bar T}$.
Hence $L$ being ample implies
that both $L_0$ and $D\cap \bar T$ are ample.

For the converse, recall that ampleness is equivalent
to being positive.
Since $\bar A\to T$ is a topologically trivial bundle,
the positivities of $c_1(L_0)$ and $c_1(D|_{\bar T})$ imply
the positivity and therefore the ampleness of $L$.
\end{proof}

\subsection{Very ampleness criterion}
We keep the notation in the previous section.

For every ample line bundle $L$ on a projective
manifold $X$ there exists a number $m$, depending on both $X$ and
the polynomial $P_L(k)=\chi(X,L^k)$ such that $L^m$ is very ample
(Matsusaka's theorem, see [Ma72]).
It is well-known that $m$ can be chosen as $3$ if $X$ is an abelian
variety (see, e.g., [Mu70]).
Here we prove a similar result for the case where $X$ is a smooth
equivariant compactification of
a semi-abelian variety.

We will employ the following auxiliary fact on line bundles
on toric varieties:

\begin{lemma}
\label{ample-on-toric}
Let $L$ be an ample line bundle on a smooth compact
toric variety $\bar T=F$ and let $p\in \bar T$.
Then there exists an effective $T$-invariant divisor $D$ with
$p\not\in \supp D$ and $L=L(D)$.
\end{lemma}
\begin{proof}
On a toric variety every ample line  bundle  is already very ample.
Thus there are sections $\tau \in H^0(\bar T,L)$ with
$\tau(p)\ne 0$.
Now consider the $T$-action on the $H^0(\bar T,L)$.
This is a linear action of the commutative reductive
group $T\cong (\C^*)^t$ and thus completely diagonalizable.
Then by the general theory of linear algebraic groups
(cf., e.g., [Bo91])
 there is a ``character'' $\chi$, i.e., a Lie group
homomorphism $\chi:T\to\C^*$ and a section $\sigma\in H^0(\bar T,L)$
such that $\sigma(p)\ne 0$ and $\sigma(t\cdot x)=\chi(t)\sigma(x)\
\forall x\in\bar T,t\in T$.
Now let $D=\div (\sigma)$, the divisor defined by $\sigma$.
\end{proof}

\begin{theorem}
\label{three-very-ample}
Let $F$ be a toric variety and let $L$ be an
ample line bundle on an equivariant compactification $\bar A$
of a semi-abelian variety $A$ of type $F$.
Then $L^{3}$ is very ample on $\bar A$.
\end{theorem}

\begin{proof}
By Propositions \ref{decomp} and  \ref{ample-decomp},
$L=\pi^*L_0\tensor L(D)$; here $L_0$ is ample on $M$
and $D$ is an $A$-invariant divisor on $\bar A$ such that
$D\cap F$ is ample on $F$ for every fiber $F$ of $\bar\pi$.

By a standard result on line bundles on abelian varieties
$L_0^3$ is very ample on $M$.

\begin{claim}
\label{ext}
Let $p\in M$ and $F=\bar\pi^{-1}(p)$. Then every section
$\sigma\in H^0(F,L^3)$ extends to a section
of $L^3$ on $\bar A$.
\end{claim}

Recall that $F\cong\bar T$ is a toric variety.
Hence $L|_F$ being ample implies that $L|_F$ is already very
ample. Furthermore $H^1(F,{\mathcal O})=\{0\}$. As a consequence
the connected Lie group $T$ acts trivially on $\Pic(F)$.
For this reason $L|_F$ is $T$-invariant and we obtain a $T$-action
on $H^0(F,L)$. Similarly we obtain a $T$-action on $H^0(F,L^3)$.
Let $\sigma\in H^0(F,L^3)$. Then there are characters
$\chi_i:T\to\C^*$ and sections $\sigma_i\in H^0(F,L^3)$
such that $\sigma=\sum_i\sigma_i$ and $t^*\sigma_i=\chi_i(t)\sigma_i$
for all $t\in T$. Now $\div(\sigma_i)$ is a $T$-invariant divisor
on $F$ for every $i$. These divisors extend to $A$-invariant divisors
$D_i=A\cdot \div(\sigma_i)$ on $\bar A$ and $\sigma_i$ extend
to sections $\tilde\sigma_i\in H^0(\bar A,L(D_i))$ with
$\div (\tilde\sigma_i)=D_i$.
Now $L(3D-D_i)$ is topologically trivial for each $i$.
Hence there are topologically trivial line bundles $H_i'\in \Pic^0(M)$
such that $L(3D-D_i)\cong\bar\pi^*H_i'$. Since $\Pic^0(M)$ is a divisible
group, there are line bundles $H_i\in \Pic^0(M)$ with $H^3_i\cong H_i'$.
The ampleness of $L_0\in \Pic(M)$ implies that $L_0\tensor H_i$ is also
ample and that consequently $(L_0\tensor H_i)^3\cong L_0^3\tensor H_i'$
is very ample.
It follows that $L_0^3\tensor H_i'$ admits a section not vanishing at $p$.
Hence there exists an element
$\zeta_i\in H^0(\bar A,\bar\pi^*(L_0^3\tensor H_i'))$
such that $\zeta_i$ has no zero on $F$.
We have
\[
\tilde\sigma_i\tensor\zeta_i\in H^0(\bar A,\bar\pi^*(L_0^3\tensor H_i')
\tensor L(D_i))\cong H^0(\bar A,L^3).
\]
Thus each $\sigma_i$ and consequently also their 
sum $\sigma=\sum_i\sigma_i$
extends to a section of $L^3$ over $\bar A$.
The claim is proved.

Using Lemma~\ref{ample-on-toric} one can show that for every point
$p\in\bar A$ there is an effective
$A$-invariant divisor $D$ (depending on $p$)
and an ample line bundle $L_0$ on $M$ such that 
$L\cong\pi^*L_0\tensor L(D)$. Thus $L^3$ can be realized as the 
tensor product of a line bundle $L(3D)$ with a section $\sigma$ not
vanishing at $p$ and a line bundle $\pi^*L_0^3$
which is the pull-back of a very ample
line bundle on $M$.

Combining these facts with Claim \ref{ext}, we can easily
verify that $L^3$ is very ample on $\bar A$.
\end{proof}

\begin{lemma}\label{numberl}
Let $\bar T$ be a toric variety. Then there exists a number
$l_0$ such that for every ample line bundle $L$ over $\bar T$
the bundle $L^{l_0}\tensor K^{-1}_{\bar T}$ is ample, too.
\end{lemma}

\begin{proof}
There are finitely many $T$-stable curves $C_1,\ldots,C_r$.
By the toric Nakai criterion ([O85], Theorem 2.18) a line bundle
$L$ on $\bar T$ is ample if and only if $\deg(L|_{C_j})>0$
for all $j\in\{1,\ldots,r\}$.
It follows that every $l_0$ with
$\deg \,(K^{-1}_{\bar T}|_{C_j})>- l_0, 1\leq j \leq r$,
has the desired property.
\end{proof}

\begin{proposition}
\label{kan-ample}
Let $A\hookrightarrow\bar A$ and $\bar T$ be as above.
Then there exists a number $l$ depending only on the toric variety
$\bar T$ such that for every ample line bundle $L$ on $\bar A$,
the line bundle $L^l\tensor K^{-1}_{\bar A}$ is ample.
\end{proposition}

\begin{proof}
By Propositions \ref{decomp} and \ref{ample-decomp}
there is an ample line bundle $L_0$ on $M$ and an $A$-invariant divisor
$D\in\Div(\bar A)$ such that $D|_F$ is ample  for every fiber $F$
of $\bar\pi$.
If $v_1,\ldots,v_m$ is a basis for $\Lie A$, then
\[
K_{\bar A}^{-1}=\div(v_1\wedge\ldots\wedge v_m).
\]
Thus $K_{\bar A}^{-1}=\partial A$.
In particular the anticanonical bundle
$K_{\bar A}^{-1}$ is induced by an $A$-invariant divisor.
The adjunction formula implies $K_{\bar A}^{-1}|_F=K_F^{-1}$.

Let $l_0$ be a number as in  Lemma~\ref{numberl}.
Then $l_0D|_F-K_{\bar T}$ is an ample divisor on the fiber $F$ of
$\bar A\to M$.
Thus $l_0D-K_{\bar A}$ is an $A$-invariant divisor on
$\bar A$ whose restriction to every fiber is ample.
On the other hand, $(L_0)^{l_0}$ is ample on $M$.
Thus
\[
(L_0)^{l_0}\tensor L(l_0D-K_{\bar A}) \cong L^{l_0}\tensor K_{\bar A}^{-1}
\]
is ample on $\bar A$ by Proposition \ref{ample-decomp}.
\end{proof}

\begin{corollary}
\label{power-ample}
For every toric variety $\bar T$ there exists a number $l$ such that
for every semi-abelian variety $A$ with smooth equivariant
compactification $A\hookrightarrow \bar A$ of type $\bar T$
and every ample line bundle $L$ on $\bar A$ the
tensor power $L^l$ is very ample and moreover
$h^0(L^k)=\chi(L^k)$, $k \geq l$.
\end{corollary}

\begin{proof}
We may choose $l$ such that $L^l$ is very ample
(Proposition~\ref{three-very-ample})
and that furthermore $L^l\tensor K_{\bar A}^{-1}$ is ample
(Proposition~\ref{kan-ample}).
By the Kodaira Vanishing Theorem, the
ampleness of $L^k\tensor K_{\bar A}^{-1}$
with $k\geq l$ implies $h^0 (L^k)=\chi(L^k)$.
\end{proof}

\section{Parametrizing spaces}
\subsection{Characterizing compactifications of semi-abelian varieties}
We start with a preparatory lemma.
\begin{lemma}\label{auto-semi-ab}
Let $A$ be a semi-abelian variety.
Then $\Aut(A)^0=A$,
where $\Aut(A)^0$ denotes the connected component of the group
of all variety automorphisms of $A$.
\end{lemma}

\begin{proof}
We have to show that every automorphism $\phi\in \Aut(A)^0$
is given as translation by some element $a\in A$.
It suffices to show that $\phi=\id_A$ if $\phi\in \Aut(A)^0$ and
$\phi(e_A)=e_A$.

Let
\[
1 \to T \to A \stackrel{\pi}{\to} M \to 1
\]
be a presentation of $A$, where
$M$ is an abelian variety and $T\cong(\C^*)^t$.
Observe that every algebraic morphism from $\C^*$ to
an abelian variety is necessarily constant.
Therefore:
If $\phi$ is an automorphism of the algebraic variety $A$,
then $\phi$ must map fibers of $\pi:A\to M$ into fibers
 and thereby
induce an automorphism $\phi_0$ of $M$.
Now $\phi\in\Aut(A)^0$ implies that $\phi$ acts trivially on the
fundamental group $\pi_1(A)$ and therefore $\phi_0$ acts trivially
on the fundamental group $\pi_1(M)$.
Together with $\phi_0(e_M)=e_M$ this implies that $\phi_0=\id_M$.
Thus $\phi$ can act only along the fibers of $\pi$.
Moreover, since the fundamental group of the fiber $T$ injects into
the fundamental group of $A$, it is clear that for each fiber
of the projection map $\pi:A\to M$ the restriction of $\phi$ is homotopic
to the identity. Now the only automorphisms of $(\C^*)^t$
which acts trivially on the fundamental group $\pi_1((\C^*)^t))\cong\Z^t$
are translations.
Hence $\phi$ must be given in the form
\[
\phi: x \mapsto \zeta(\pi(x))\cdot x,
\]
where $\zeta:M \to T$ is a holomorphic map.
But $M$ is compact, hence every holomorphic map from $M$ to
$(\C^*)^t$ is constant. Thus $\zeta$ is constant, and in fact
$\zeta\equiv e_A$, because $\phi(e_A)=a_A$, i.e., $\phi=\id_A$.
\end{proof}

\begin{proposition}\label{characterize-semi-ab}
Let $X$ be a smooth projective variety with $n=\dim X$
and let $D$ be an effective reduced divisor on $X$.
Let $V$ denote the vector space of those vector fields on $X$
which are everywhere tangent to $D$.
Then $X$ is a equivariant algebraic compactification of
a semi-abelian variety $A$ with $A=X\setminus D$
if and only if the
following conditions are fulfilled:
\begin{enumerate}
\item
$\dim V=n=\dim X$.
\item
For every $v,w\in V$ the vector fields $v$ and $w$ commute.
\item
Let $v_i, 1 \leq i \leq n$, be bases of $V$.
Then the zero divisor of
$\wedge_{i=1}^n v_i$ equals $D$.
\item
For $v \in V$
and $p\in \{x\in X:v_x=0\}$ let $\Phi_{v,p}$ denote the
induced linear endomorphism of $T_pX$.
Then $\Phi_{v,p}$ is nilpotent if and only if $v=0$.
\end{enumerate}
\end{proposition}
\begin{proof}
Conditions (ii) and (iii) ensure that $\Aut(X)$ contains a connected
commutative subgroup $A_0$ (corresponding to the Lie algebra $V$)
which has an open orbit $\Omega$, and that this
open orbit $\Omega$ is precisely $X\setminus D$.
The connected component of $\Aut(X)$ is an algebraic group.
Let $A$ denote the Zariski closure of $A_0$ in $\Aut(X)$.
Since $A_0$ stabilizes $D$, its Zariski closure $A$ must stabilize
$D$, too.
But this implies $\Lie(A_0)\subset V$.
Hence $A=A_0$.

For $p\in\Omega$ let $A_p=\{a\in A:a(p)=p\}$.
Because $A$ is commutative, $A_p=A_q$ for all $p,q\in\Omega$.
Since $\Omega$ is dense in $X$, it follows
that $A_p=\{\id_X\}$.
Hence $\Omega\cong A$.

As a connected commutative algebraic group, either $A$ is a semi-abelian
variety, or it contains an algebraic subgroup $U$ isomorphic to
the additive group $(\C,+)$. If $A$ contains such a subgroup $U$,
then every non-trivial $U$-orbit is affine and one-dimensional,
and therefore contains a $U$-fixed point $p\in X$ in its closure.
Being unipotent,
elements of $U$ induce  unipotent endomorphisms of $T_pX$.
This however implies that a vector field $v\in \Lie(U)$
induces a nilpotent endomorphism of $T_pX$, contradicting
property (iv).
Thus $A$ cannot contain an algebraic subgroup isomorphic to
the additive group and therefore must be a semi-abelian
variety.

Conversely, let $A\hookrightarrow\bar A$
be a smooth algebraic equivariant
compactification of a semi-abelian variety $A$.
Let $G$ denote the group of all automorphisms of $\bar A$
stabilizing the boundary $\partial A$.
Then every $g\in G$ stabilizes $A$ as well, hence $G\hookrightarrow
\Aut(A)$.
But $\Aut(A)^0=A$ (Lemma \ref{auto-semi-ab}).
Hence $G^0=A$, and the fundamental vector fields
of the $A$-action on $\bar A$ are the only vector fields
on $\bar A$ which are everywhere tangent to the divisor
$D=\partial A$.
This implies properties (i)--(iii).
Property (iv) follows from the fact that all the isotropy groups
of the $A$-action on $\bar A$ are algebraic subgroups and therefore
reductive.
\end{proof}

Thus, we obtain the following:
\begin{theorem}
Let $p:U\to S$ be a flat family of smooth projective varieties,
let $\Theta$ be an effective divisor on $U$,
and let $S_0$ be the set of points $s \in S$ such that
$p^{-1}(s)\setminus\Theta$
is a semi-abelian variety with $p^{-1}(s)\setminus\Theta\hookrightarrow
 p^{-1}(s)$
as a smooth equivariant compactification.
Then $S_0$ is a Zariski open subset of its Zariski closure
in $S$.
\end{theorem}

\begin{proof}
By the preceding Proposition \ref{characterize-semi-ab},
the condition $s\in S_0$ translates
into four conditions (i), $\ldots$, (iv),
each of which is closed,
or at least locally closed for a flat family.
Hence the assertion follows.
\end{proof}

\subsection{Curves}
Let $C$ be a smooth compact complex algebraic curve of genus $g$,
and let $E$ be a divisor on $C$ with $\deg (E)=2g+1$.
Then $\deg (K_C-E)=-3$, since $\deg (K_C)=2g-2$.
By Riemann-Roch's theorem, $E$ is very ample,
\[
h^0(C, L(E))=\dim H^0(C, L(E))=\deg(E) +1-g=g+2,
\]
and
\[
\chi(C, L(lE))= l\deg(E)+1-g= l(2g+1)+(1-g)
\]
for all $l\in\N$.
Hence every smooth compact algebraic curve of genus $g$ can be embedded
into $\P^{g+1}$ with Hilbert polynomial $P_1(l)=l(2g+1)+(1-g)$.

\begin{proposition}
\label{U_1}
Let $g,d\in\N$.
Then there exist projective algebraic varieties
$U_1$ and $Q_1$, a proper surjective morphism $\pi_1:U_1\to Q_1$,
and a divisor $\Sigma_1\subset U_1$ such that
for every algebraic curve $C'$ of genus $g$ with $d$ punctures,
there exists a point $p\in Q_1$ (not necessarily unique)
such that $(\pi_1^{-1}(p)\setminus \Sigma_1) \cong C'$.
\end{proposition}
\begin{proof}
Let $\pi: U\to Q_0$ be the universal family of the Hilbert scheme
of compact algebraic curves of genus $g$
with the Hilbert polynomial $P_1(l)$ as above.
Define
\[
\pi_1: U_1 =\{(u_0, u_1,\ldots,u_d)\in U^{d+1}:
\pi(u_0)=\ldots=\pi(u_d)\}\to
Q_0\times U^d=Q_1
\]
by $\pi_1(u_0,\ldots,u_d)=(\pi(u_0),u_1,\ldots,u_d)$.
Moreover we define $\Sigma_1$ as the {\em ``diagonal divisor''}
by
\[
\Sigma_1=\left\{(u_0, u_1,\ldots,u_d)\in U_1: u_0\in \{u_i\}_{i=1}^d
\right\}.
\]
\end{proof}

\subsection{Semi-abelian varieties}
\begin{proposition}
\label{U_2}
Let $(\C^*)^t\cong T\hookrightarrow\bar T$
be a smooth projective toric variety, let $\{\alpha_i\}_{i\in I}$
be a $\Z$-module basis of $ H_{2*}(\bar T,\Z)$
with $\alpha_i\in H_{2d_i}(\bar T,\Z)$ ($d_i\in\{0,\ldots,\dim T\}$)
and let $s_i\in \N$ for $i\in I$.
Then there exist projective algebraic varieties $U_2$ and $Q_2$,
a proper surjective morphism $\pi_2:U_2\to Q_2$, divisors
$\Sigma_2,\Theta\in \Div(U_2)$
such that for every semi-abelian variety $A$ of dimension $n$ ($\geq t$)
with smooth equivariant compactification $A\hookrightarrow\bar A$
of type $\bar T$ we have the following property:
If there exists an ample divisor $D\in\Div(\bar A)$
with $\alpha_i\cdot[D]^{n-d_i}=s_i$ for all $i\in I$,
then there exists a point $q\in Q_2$ and an isomorphism
$\zeta:\bar A\to \pi_2^{-1}(q)\subset U_2$ such that
$A=\bar A\setminus \zeta^{-1}(\Sigma_2)$
and $D=\zeta^{-1}(\Theta)$.
\end{proposition}

\begin{proof}
First note that topologically $\bar A\cong  \bar T \times M$ with
$M=A/T$
and that $H_*(\bar T,\Z)$ is generated by the invariant divisors
of $\bar T$ while $T_M$ is trivial.
Therefore the Chern classes of $\mathbf{T}_{\bar A}$ can be expressed
in terms of the generators of $H_2(\bar T, \Z)$.
Via Hirzebruch-Riemann-Roch's theorem
it follows that the Hilbert polynomial
of $D$ is determined by the conditions
$\alpha_i\cdot[D]^{n-d_i}=s_i$ ($i\in I$).
Furthermore $\bar T$ determines a number $l$ such that
$lD$ is very ample with $h^0(L(lD))=\chi(L(lD))$
(Corollary~\ref{power-ample}).
Hence these conditions also determine a polynomial $P_2$ and a number $k$
such that $\bar A$ can be embedded into $\P^k(\C)$
with Hilbert polynomial $P_2$.

The boundary $\partial A$ of $A$ in $\bar A$ is a divisor and its
intersection numbers are determined by the choice of the toric
variety $\bar T$. Thus the Hilbert polynomial of this divisor
$\partial A$ is also determined by the choice of the toric variety.

By the general theory of Hilbert schemes we now obtain a morphism
between projective algebraic varieties $\pi_2:U_2\to Q_2$
and divisor $\Sigma_2, \Theta\in \Div(U_2)$ satisfying the following
 property:
\begin{description}
\item[  ] \hskip10pt
For every $(A,\bar A,D)$ with the given numerical data we can find
a point $q\in Q_2$ such that there is
an isomorphism $\zeta : \bar A \to  \pi_2^{-1}(q)$
with $\zeta^{-1}(\Theta)=D$ and $\zeta^{-1}(\Sigma_2)=\partial A$.
\end{description}

By Proposition \ref{characterize-semi-ab}
we may replace $Q_2$ by a closed subspace and therefore
assume that there is a Zariski-open subset $W\subset Q_2$
such that the fiber $\pi_2^{-1}(q)$ is a smooth equivariantly
compactified semi-abelian variety for every $q\in W$.
\end{proof}

\section{Proof of the Main Theorem in general case}

\subsection{Compactified total space of a coherent sheaf}\label{grauert}

If $E\to X$ is a vector bundle of rank $n$ over a compact complex space,
then $E$ admits a compactification in the form of an
embedding into a $\P^n$-bundle $\bar E$ over $X$.
This construction can be generalized to coherent sheaves.

Following Grauert [G62], we define a {\em linear space} over $X$ as
a holomorphic map of complex spaces $\pi: V\to X$ together
with holomorphic maps $\mu :V \times_X V  \to V$ and
$\nu:\C\times V\to V$, satisfying the conditions:
\begin{enumerate}
\item
$\pi(\mu(x,y))=\pi(x)=\pi(y)$ for all $(x,y)\in V\times_XV $
and $\pi(\nu(t,v))=\pi(v)$ for all $t\in\C$ and $v\in V$.
\item
Each fiber of $\pi$ is a complex vector space, where the vector
addition is given by $\mu$ and
the scalar multiplication by $\nu$.
\item
For every point $x\in X$ there is an open neighbourhood $U$
and a number $n\in\N$ such that
there is a commutative diagram
\[
\begin{matrix}
V|_U & \stackrel{\phi}{\to} & \C^n \times U \\
\quad\downarrow {}_\pi & &\quad\quad \downarrow \mathrm{}_{\mathrm{pr}_2}\\
U & \stackrel{\id}{\to} & U,
\end{matrix}
\]
where $V|_U=\pi^{-1}(U)$ and $\phi$ is an embedding,
linear on every fiber
of $\pi$.
\end{enumerate}

As described in [G62],
for every coherent sheaf $\mathcal{S}$ on a complex space $X$
there exists a ``{\em linear space}'' $V\to X$ such that
$\mathcal S$ is isomorphic to the sheaf of those holomorphic
functions on $V$ which are linear on each fiber.
Furthermore, if $X$ is compact, there is a natural compactification
$V\hookrightarrow\bar V$
induced by $\C^n\hookrightarrow\P^{n}$ in each fiber.
The points in the boundary $\partial V=\bar V\setminus V$
correspond to complex lines in the fibers of $\pi$.

\subsection{Proof of the Main Theorem}

Let $U_1\stackrel{\pi_1}{\to} Q_1$ and $U_2\stackrel{\pi_2}{\to} Q_2$,
together with $\Sigma_i\in \Div(Q_i)$ and
$\Theta\in\Div(Q_2)$ be as constructed
in Propositions \ref{U_1} and \ref{U_2}, respectively.
We set
\begin{align*}
\pi=\pi_1\times\pi_2 &: U=U_1\times U_2 \to Q=Q_1 \times Q_2,\\
p_i &: Q \to Q_i,
\end{align*}
where $p_i$ are the natural projections.
We define $Q_i^*$ as the open subset
of $Q_i$ consisting of all those points $q_i\in Q_i$
such that
\begin{enumerate}
\item
the fiber $\pi_i^{-1}(q_i)$ is irreducible,
reduced and smooth,
\item
$\pi_i^{-1}(q_i)\not\subset \Sigma_i$, $i=1,2$,
and $\pi_2^{-1}(q_2)\not\subset \Theta$.
\end{enumerate}
Set
$$
Q^*=Q_1^* \times Q_2^*.
$$
The points of $Q^*$ are the \begin{em}``nice''\end{em} points
in which we are really interested.

For each of the morphisms $\pi_i:U_i \to Q_i$,
 the sheaf of horizontal logarithmic differential
forms along $\Sigma_i$ is defined to be
the subsheaf of those logarithmic differential forms whose induced
differentials on fibers vanish.
This is a coherent sheaf as well as the
quotient sheaf of all differential forms divided by the horizontal
ones.
This quotient sheaf we call the ``sheaf of vertical
logarithmic differential forms''
denoted by $\Omega^1_{\mathrm{vert}} (U_i,\log \Sigma_i)/Q_i$.
Define
$$
\sS_i=p_i^{-1}{\pi_i}_*\Omega^1_{\mathrm{vert}}
 (U_i,\log \Sigma_i)/Q_i.
$$
This is a coherent sheaf over $Q$.
Moreover, it is a locally free sheaf over $Q^*$.
By Grauert's theory of linear spaces (see \S\ref{grauert}),
there are linear
spaces $S_i\to Q$ such that $\sS_i$ is the sheaf of functions on $S_i$
which are linear on fibers.

Set
$$
S_i^*=S_i|_{Q^*},\quad U_i^*=U_i|_{Q^*}, \quad i=1,2.
$$
The fibers of $S^*_i\to Q^*$ can be identified with
the space of logarithmic vector fields on the quasi-Albanese variety
of the fibers of $U^*_i \to Q^*_i$.
Let $J^k_{\mathrm{vert}}(U_i^*)$ denote the relative
$k$-jet space, i.e., the inverse image of zero jets of
$\pi_{i_*}:J^k (U^*_i)\to J^k(Q^*_i)$.
Let $L\to Q$ be the space of all relative linear morphisms from
$S_1$ into $S_2$ over $Q$.
We set $L^*=L|_{Q^*}$.

By Proposition \ref{jet-trivial} we can identify
$\m(k)\tensor S_{iq_i}^*$
with $J^k(U_{iq_i}^*,\log\Sigma_{iq_i})$ for $q_i \in Q^*_i$,
and
$$
J^k(U_{iq_i}^*,\log\Sigma_{iq_i})|_{U^*_{iq_i}\setminus \Sigma_{iq_i}}
\cong
J^k_{\mathrm{vert}}(U_i^*)|_{U^*_{iq_i}\setminus \Sigma_{iq_i}}.
$$
As in (\ref{k-locus}), we next define algebraic subsets
$Z_k \subset L^*\times_{Q^*}( (U^*_1\setminus \Sigma_1)
\times  (\Theta^*\setminus \Sigma_2))$, $k=1,2, \ldots$ by
\begin{align*}
Z_k=\Big\{ & (\lambda,p,q)\in
 L^*\times_{Q^*}( (U^*_1 \setminus \Sigma_1)
\times  (\Theta^*\setminus \Sigma_2)): \\
& \left(\mathrm{id}_{\m(k)}\tensor\lambda\right) \left(
J^k_{\mathrm{vert}}(U^*_1)_{p_1(p)}\right)
\subset J^k_{\mathrm{vert}}(\Theta^*)_{p_2(q)}
\Big\}.
\end{align*}
Since $\{Z_k\}_k$ is a decreasing sequence of algebraic subsets,
it terminates.
Thus there exists a number $N$ such that $Z_k=Z_N$
for all $k\ge N$.

Let $q\in Q^*$, $C=\pi_1^{-1}(p_1(q))\setminus\Sigma_1$,
$A=\pi^{-1}(p_2(q))\setminus\Sigma_2$ and
$D=A\cap\Theta$.
Let $f:C\to A$ be a morphism.
Then $C$ may be regarded as a closed subvariety of $\Alb(C)$,
and there is a morphism $\tilde f: \Alb(C) \to A$ such
that $\tilde f|_C=f$.
If $\mult_x\, f^*D \geq k$ with $x\in C$, then by construction
$(d\tilde f_x, x, f(x)) \in Z_k$.
If $k \geq N$, then $f(C)$ is tangent to $D$ at $f(x)$ with infinite
order; that is, $f(C) \subset D$.
This finishes the proof of the Main Theorem.
{\it Q.E.D.}\par\vskip3pt
\section{Application}

As an application of the Main Theorem we prove

\begin{theorem}
\label{appl}
Let $C$ be an affine
algebraic curve, let $A$ be an abelian variety,
and let $D$ be an ample  effective reduced divisor on $A$.
Then either there is a non-constant morphism from $C$ into $D$
or there are only finitely many morphisms from $C$ into
$A\setminus D$.
\end{theorem}

Before we prove the theorem we need the following.

\begin{lemma}\label{local-rigid}
Let $C$ be a smooth algebraic curve with
smooth algebraic compactification
$C\hookrightarrow \bar C$,
let $A$ be an abelian variety,
and let $D$ be an ample hypersurface in $A$.
Let $f:\bar C\to A$ be a morphism with $f(C)\subset A\setminus D$,
and let $f_a:C\to A$ be given by $f_a(t)=a+f(t)$.
Set
\[
P=\{a\in A: f_a(C)\subset A\setminus D\}.
\]
Then either there exists a point $p\in A$ with $f_p(C)\subset D$
or $P$ is finite.
\end{lemma}
\begin{proof}
Let $S=\bar C\setminus C$.
Let $d=\deg f^*D$. Let $(S_i)_{i\in I}$ denote the finite family
of all distinct effective divisors on $\bar C$ with $\deg(S_i)=d$
and $\supp S_i \subset S$.
Define
\[
P_i=\{a\in A: f_a^*(D)=S_i\}
\]
for $i\in I$,
and
\[
P_\infty=\{a\in A:f_a(C)\subset D\}.
\]
Then $P$ is the disjoint union of $P_\infty$
and $P_i, i \in I$.
Assume that $P_\infty=\emptyset$.
Then $P_i$ is compact.
Now consider
\[
F_i: P_i \times C\to A\setminus D
\]
given by $F_i(a,t)=f_a(t)=a+f(t)$.
Note that $A\setminus D$ is affine, because $D$ is assumed to be ample.
Thus the compactness of $P_i$ implies that
$F_i$ is locally constant with respect to the first variable.
Due to the equality $F_i(a,t)=a+f(t)$ this can be true only
if $P_i$ is finite.
Hence $P$ is finite.
\end{proof}

\begin{corollary}\label{cor.8.3}
Assume that every morphism from $C$ to $D$ is constant.
Then $\Mor(C,A\setminus D)$ is embedded into $\Mor(\bar C,A)$
as a discrete subset.
\end{corollary}
\begin{proof}
Let $f:C\to A\setminus D$ be a morphism, and let $U$ be a Stein open
neighbourhood of the neutral element $e$ in $A$.
Then $W=\{g:\bar C\to A: f(x)-g(x)\in U,\: \forall x\in\bar C\}$
is an open neighbourhood (for the compact-open topology)
of $f$ in $\Mor(\bar C,A)$.
However, if $f(x)-g(x)\in U$ for all $x\in\bar C$, then
$x\mapsto f(x)-g(x)$ is constant, because $\bar C$ is compact and
$U$ is Stein.
Thus every $g\in W$ is given as $g:x\mapsto a+f(x)$ for some $a\in A$.
Therefore the assertion follows from above Lemma \ref{local-rigid}.
\end{proof}

{\it Proof of Theorem \ref{appl}.}

Assume that every morphism from $C$ to $D$ is constant.
By the Main Theorem this implies that there exists a number $N$
such that for every morphism $\bar f:\bar C\to A$ the multiplicities of
$\bar f^*D$
are uniformly bounded by $N$.
This leads to $\deg\bar f^*D\le N\cdot\#(\bar C\setminus C)$
for all holomorphic maps $\bar f:\bar C\to A$ with
$\bar f(C)\subset A\setminus D$.
Thus $\Mor(C, A \setminus D)$ is contained in a union of finitely
many irreducible components of $\Mor(\bar C,A)$.
It follows from Corollary~\ref{cor.8.3} that
$\dim \Mor(C, A \setminus D)=0$. Therefore $\Mor(C, A \setminus D)$ must
be finite.
{\it Q.E.D.}

\begin{remark}{\rm
\begin{enumerate}
\item
By Dethloff-Lu [DL01],
every holomorphic mapping $f:C \to A\setminus D$ extends
holomorphically to $\bar f:\bar C \to A$,
and so $f$ is an algebraic morphism.
Thus, Theorem \ref{appl} implies that there are only finitely
many non-constant holomorphic mappings of $C$ into $A\setminus D$.
\item
If $D$ contains no translate of a non-trivial Abelian subvariety
of $A$, then $A\setminus D$ is complete hyperbolic and hyperbolically
embedded into $A$ ([Gr77]).
In this case, the uniform bound for
$\deg \bar f, f \in \Mor(C, A \setminus D)$,
is a consequence of [Nog88].
\item
If $D$ admits a non-constant morphism from $C$, then there may be
infinitely many non-constant morphisms $f: C \to A$
with $f^{-1}D \subset S$.
In fact, let $E$ be an elliptic curve, and set $A=E^2$.
Then $D=\{0\}\times E+E\times\{0\}$ is ample.
Set $\bar C=E$ and $C=E \setminus \{0\}$.
There are infinitely many morphisms
$$
f_a: x \in C \to (x, a)\in A\setminus D,
\qquad a \in C.
$$
In this case, $\bar f_0(C) \subset D$.
\end{enumerate}}
\end{remark}

\bigskip

\rightline{Noguchi, Junjiro}
\rightline{Graduate School of Mathematical Sciences}
\rightline{University of Tokyo}
\rightline{Komaba, Meguro,Tokyo 153-8914}
\rightline{e-mail: noguchi@ms.u-tokyo.ac.jp}
\medskip

\rightline{Winkelmann, J\"org}
\rightline{Korea Institute for Advanced Study}
\rightline{School of Mathematics}
\rightline{207-43 Cheongryangri-dong Dongdaemun-gu}
\rightline{Seoul 130-012, Korea}
\rightline{ e-mail: jwinkel@member.ams.org}
\rightline{Web: http://www.math.unibas.ch/\~{ }winkel/index.html}

\end{document}